\documentclass[12pt]{article}
\usepackage{amssymb}
\def \C{{\mathbb C}}

\def \Z{{\mathbb Z}}
\def \p{{\mathbb P}}

\begin{document}
\title{Vari\'et\'es de Fano et morphismes: notes d'un 
mini-cours}
\author{Ekaterina Amerik}

\date{}
\maketitle

Ce texte est une version partielle d'un mini-cours que j'ai fait
lors d'une rencontre de l'ACI Jeunes Chercheurs ``Dynamique des
applicatons polynomiales'' \`a Toulouse en novembre 2004. Le mini-cours,
tel que je l'ai prepar\'e, 
\'etait en trois parties, dont deux sont present\'ees ici. 
Dans la premi\`ere, on expose quelques pr\'eliminaires
sur les vari\'et\'es de Fano de dimension trois. Tous les faits qu'on discute ici
peuvent se trouver dans la litterature, et beaucoup dans \cite{I}; mais 
d\'ej\`a le travaux \cite{I} sont assez volumineux, et leur but est tout \`a fait
diff\'erent du n\^otre (qui, lui, est bien moins ambitieux). De plus, pour 
quelques exemples classiques, il est plut\^ot difficile d'indiquer la meilleure
r\'ef\'erence.

La seconde partie pr\'esente quelques m\'ethodes pour borner le degr\'e
d'un morphisme $f:X\rightarrow Y$, o\`u $X, Y$ sont des vari\'et\'es lisses
projectives de dimension 3, avec le groupe de Picard cyclique, et $Y$
est une vari\'et\'e de Fano diff\'erente de $\p^3$. Les m\'ethodes sont
extraites de quelques articles assez recents (\cite{A1},\cite{A2},\cite{ARV}).

La troisieme partie n'est pas reprise dans ces notes, car pour cette partie-ci,
il existe une r\'ef\'erence pr\'ecise, assez courte et {\it ``self-contained''}:
l'article \cite{HM}. Dans le mini-cours, on a essay\'e d'expliquer le th\'eor\`eme suivant
de Hwang et Mok: {\it soient $X, Y$ des vari\'et\'es de Fano du groupe de 
Picard cyclique, dimension quelconque. Si $Y$ contient des courbes au fibr\'e
normal trivial, alors le degr\'e de $f:X\rightarrow Y$ est born\'e}.

Voici quelques notations utilis\'ees:

- Pour $X$ une vari\'et\'e telle que $Pic(X)=\Z$, $H_X$ (ou simplement $H$
s'il n'y a pas de risques de confusion) d\'enote le g\'en\'erateur ample
du $Pic(X)$;

- $K_X$ d\'enote le diviseur canonique, ou, par abus de notation, le
fibr\'e en droites canonique;

- $h^i(X,E)=dim(H^i(X,E))$;

-$|H|$ le syst\`eme lin\'eaire associ\'e au fibr\'e en droites $H$;

- $\Omega_X$ le fibr\'e cotangent de $X$

\

{\bf Partie I: Preliminaires sur les vari\'et\'es de Fano de dimension 3}

\

{\bf Definition:} Une vari\'et\'e lisse projective $X$ est dite {\it de
Fano}, si le fibr\'e anticanonique ${\cal O}(-K_X)=
det(T_X)$ est ample. On appelle {\it l'indice} d'une
vari\'et\'e de Fano $X$ le plus grand entier $r$ tel que
${\cal O}(-K_X)=L^{\otimes r}$ pour un fibr\'e en droites
$L$ sur $X$.

\

L'indice d'une vari\'et\'e de Fano ne peut pas \^etre
tr\`es grand, selon le th\'eor\`eme suivant, d\^u
\`a Kobayashi et Ochiai (\cite{KO}):

\

{\bf Th\'eor\`eme:} {\it L'indice $r$ d'une vari\'et\'e de
Fano $X$ de dimension $n$ n'excede pas $n+1$. De plus,
si $r=n+1$, alors $X\cong \p^n$; si $r=n$, alors
$X$ est isomorphe a une quadrique dans $\p^{n+1}$.}

\

Dans ce mini-cours, on va s'int\'eresser aux vari\'et\'es
de Fano dont le groupe de Picard est $\Z$. Notons que
ceci est la m\^eme condition que $b_2=1$: en effet,
par le th\'eor\`eme d'annulation de Kodaira, 

$$H^i(X,{\cal O}_X)=0\ pour \ i>0,$$
on a donc par la th\'eorie de Hodge $H^2(X,\C)=H^{1,1}(X)$,
et l'on conclut par le th\'eor\`eme de (1,1)-classes de Lefschetz. 

Le premier cas int\'eressant est $n=dim(X)=3$. Effectivement,
il est \'evident que la seule vari\'et\'e de Fano de dimension un
est $\p^1$; quant a la dimension deux, il n'est pas difficile
\`a d\'emontrer la proposition suivante:

\

{\bf Proposition:} {\it Si $X$ est Fano, $dim(X)=2$ et $Pic(X)=\Z$, alors $X\cong \p^2$}.

\

Effectivement, $c_2(X)=\chi(X)=4$ par l'annulation de nombres de Hodge
ci-dessus, et $c_1^2(X)+c_2(X)=12\chi({\cal O}_X)=12$ par la formule
de Riemann-Roch $\chi(X,E)=(ch(E)\cdot td(X))_n$: en dimension deux, cette 
formule s'\'ecrit de la fa\c{c}on
suivante: $$\chi({\cal O}_X(D))=D(D-K)/2+(K^2+c_2)/12$$ pour $D$ un
diviseur sur $X$. D'o\`u $K^2=9$. Comme la forme d'intersection
sur $H^2(X,\Z)=Pic(X)$ est unimodulaire (dualit\'e de Poincar\'e),
on voit que l'indice de $X$ est 3: $K=-3H$, o\`u $H$ est le g\'en\'erateur
ample de $Pic(X)$. M\'eme si l'on oublie le th\'eor\`eme de Kobayashi-Ochiai
(dont la preuve est par ailleurs une g\'en\'eralisation, bien que non-\'evidente, 
du petit argument
qui suit!),
on en d\'eduit tout de suite que $X\cong \p^2$: Riemann-Roch donne
$\chi({\cal O}_X(H))=3$, l'annulation de Kodaira implique 
$$\chi({\cal O}_X(H))=dim(H^0(X,{\cal O}_X(H)),$$ $H$ d\'efinit donc
une application rationnelle $\phi: X\rightarrow \p^2$, et la condition
$H^2=1$ implique facilement que $\phi$ est un isomorphisme. Effectivement,
on doit juste voir que $|H|$ n'a pas de points de base; il n'y a clairement pas
de composantes de base parce que $H$ engendre $Pic(X)$; s'il y a un point de base,
il n'y a qu'un seul, du fait que $H^2=1$, et chaque membre de $|H|$ coupe les autres
membres seulement en ce point de base. Ceci est impossible car il y a ``trop de membres'':
$dim|H|=2$. 

\

En dimension trois, il devient beaucoup plus difficile de classifier les
vari\'et\'es de Fano dont le groupe de Picard est cyclique; ceci a \'et\'e
fait par V.A. Iskovskih (\cite{I}) \`a la fin des ann\'ees 1970. Nous ne pouvons
donner ici l'id\'ee de la preuve de cette classification, ni m\^eme la liste tr\`es
pr\'ecise des vari\'et\'es obtenues; voici sa version un peu
abreg\'ee.

\

{\it 1.1. Classification d'Iskovskih:}

\

On sait d\'ej\`a que $ind(X)=4$ implique que $X=\p^3$, et que $ind(X)=3$
implique que $X \cong Q_3$ (la quadrique). 

\

{\it Le cas ind(X)=2:} Il y a en plus cinq familles
de vari\'et\'es de Fano d'indice deux: si l'on d\'enote $H$ le g\'en\'erateur
ample de $Pic(X)$, alors $H^3$ (le ``degr\'e'' de $X$)
peut varier de 1 \`a 5. On calcule par 
Riemann-Roch 
$$\chi(X,D)=D^3/6-KD^2/4+D(K^2+c_2)/12+c_1c_2/24$$
et par l'annulation de Kodaira que $H^0(X,{\cal O}(H))=H^3+2$.
Si $d=H^3>2$, $H$ est tr\`es ample et notre vari\'et\'e se trouve donc
plong\'ee dans $\p^{d+1}$. Plus pr\'ecisement, 

- $H^3=3$: $X$ est une cubique dans $\p^4$;

- $H^3=4$: $X$ est l'intersection de deux quadriques dans $\p^5$;

- $H^3=5$: $X$ est l'intersection de $G(1,4)\subset\p^9$ et trois hyperplans
(ici, $G(1,4)\subset\p^9$ est la grassmannienne des droites dans $\p^4$,
plong\'ee dans $\p^9$ par les coordonn\'ees de Pl\"ucker). Dans ce cas, il
n'y a qu'un seul $X$, \`a un isomorphisme pr\`es: toutes les vari\'et\'es lisses
obtenues comme ci-dessus sont isomorphes.

- si $H^3=2$, $H$ est toujours sans points de base et $X$ est un rev\^etement
double de $\p^3$ ramifi\'e le long d'une quartique.

- si $H^3=1$, alors $|H|$ a un point de base et $X$ est un rev\^etement
double d'un c\^one sur la surface de Veron\'ese (=$\p^2$ plong\'e dans $\p^5$
par les formes quadratiques), ramifi\'e le long d'une section hypercubique
de ce c\^one.
 
On va parfois noter $A_k$ une vari\'et\'e de Fano d'indice 2, Picard $\Z$, 
dont $H^3=k$.

\

{\it le cas $ind(X)=1$:}

Par la formule de Riemann-Roch, on trouve $\chi(X, -K_X)=(-K_X)^3/2+3$, en
particulier, $(-K_X)^3$ est pair. On \'ecrit souvent $(-K_X)^3=2g-2$, et l'on appelle
ce $g$ le genre de $X$ (si $H=-K_X$ est tr\`es ample, c'est, bien s\^ur, le genre
de la courbe section lineaire de $X$ plong\'e anticanoniquement). On a alors $dim|H|=g+1$.
La classification affirme alors que $H^3$ prend les valeurs paires de 2 \`a 22,
except\'e 20. On d\'enote $V_k$ une vari\'et\'e de Fano d'indice 1 dont $H^3=k$.

On a les cas suivants:

- $V_2$ est un rev\^etement double de $\p^2$ (ramifi\'e le long d'une sextique);
 
- $V_4$ est ou bien une quartique de $\p^4$, ou bien un rev\^etement double
d'une quadrique dans $\p^4$ (ramifi\'e le long d'une section hyperquartique);

Si $k>4$, $-K_{V_k}$ est toujours tr\`es ample;

- $V_6$ est une intersection d'une quadrique et une cubique dans $\p^5$;

- $V_8$ est une intersection de trois quadriques dans $\p^6$;

Les $V_k$ avec $k$ plus grand ne sont pas des intersections compl\`etes, mais
on peut toujours les d\'ecrire en termes de certains fibr\'es vectoriels sur
certaines vari\'et\'es homog\`enes. Par exemple:

- $V_{10}$ est une intersection du c\^one $K$ sur $G(1,4)$ ($K\subset\p^{10}$) avec
trois hyperplans et une quadrique;

- $V_{14}$ est une intersection de $G(1,5)\subset \p^{14}$ avec 5 hyperplans;

- $V_{22}$ est le lieu des z\'eros d'une section du fibr\'e $\Lambda^2U^*\oplus
\Lambda^2U^*\oplus\Lambda^2U^*$ sur $G(2,6)$ (\cite{M}).

On d\'emontre aussi que les vari\'et\'es $V_k, k\geq 8$ sont intersections de
quadriques qui les contiennent (dans le plongement canonique).

\

{\it 1.2. Une approche aux vari\'et\'es de Fano de petite degr\'e: espaces projectifs
\`a poids}

\

{\bf D\'efinition:} L'espace projectif \`a poids 
(o\`u ``l'espace projectif anisotrope''
$\p(a_0, \dots , a_n)$ est le quotient de $\C^{n+1}-0$ par l'action suivante
de $\C^*$: $$\lambda (x_0,\dots ,x_n)=(\lambda^{a_0}x_0,\dots ,\lambda^{a_n}x_n).$$ 
(D'une fa\c{c}on plus alg\'ebrique, on la d\'efinit comme $Proj$ de l'anneau
de polyn\^omes $\C[x_0,\dots , x_n]$ gradu\'e par $deg(x_i)=a_i$.)

\

Il est clair que $\p(a_0, \dots , a_n)\cong \p(Na_0, \dots , Na_n)$ pour tout entier
$N$, on peut donc supposer que $pgcd(a_i)=1$.

On peut aussi voir $\p(a_0, \dots , a_n)$ comme quotient de $\p^n$ par l'action
\'evidente de $\mu_{a_0}\times \dots \times \mu_{a_n}$:

$$(\epsilon_0,\dots, \epsilon_n)(x_0:\dots :x_n)=(\epsilon_0x_0:\dots :\epsilon_nx_n).$$

Effectivement, l'application $\C^n\rightarrow \C^n$ donn\'ee par $x_i=y_i^{a_i}$
descend aux quotients correspondants.

\

{\it Exemple:} $\p(1,1,2)$ est le c\^one quadratique dans $\p^3$;
le plongement dans $\p^3$ est donn\'e par 
$(x_0:x_1:x_2)\rightarrow (x_0^2:x_0x_1:x_1^2:x_2)$. De m\^eme, $\p(1,1,n)$ est le c\^one
sur la courbe rationnelle normale de degr\'e $n$.

\

On peut m\^eme supposer que pour tout $i$, $pgcd(a_0,\dots,\hat{a_i}, 
\dots, a_n)=1$: sinon on ``divise'' par ce pgcd. En effet,
considerons le morphisme du quotient $\pi:\p^n\rightarrow \p(a_0, \dots , a_n)$.
Si $c=pgcd(a_1,\dots, a_n)$, on peut factoriser $\pi$ par 
$\pi_1: \p^n\rightarrow \p(1, c, \dots ,c)$, morphisme du quotient
par l'action de $\{e\}\times\mu_c^{\times n}$. Mais ceci n'est autre
chose que le morphisme de $\p^n$ vers lui-m\^eme qui monte toutes les
coordonn\'ees \`a la puissance $c$.

Il y a une maniere commode de recouvrir $\p(a_0, \dots , a_n)$
par des quotients de l'espace affine: consid\'erons, dans $\C^n-0$, 
le sous-ensemble $A_i=\{x|x_i=1\}$, stabilis\'e par le groupe $\mu_{a_i}$;
le quotient de $A_i$ par l'action de ce groupe est l'ouvert $U_i=\{x|x_i\neq 0\}$
de $\p(a_0, \dots , a_n)$. En particulier il est lisse si $a_i=1$.
On voit alors que les singularit\'es de $\p(a_0, \dots , a_n)$
sont donn\'ees par l'annulation de certaines coordonn\'ees:
$$Sing(\p(a_0, \dots , a_n))=\bigcup_{k>1}\{x|x_i=0\ si\ k\ ne\ divise\ pas\ a_i\}.$$

Soit $U(a_0, \dots , a_n)$ la partie non-singuli\`ere de 
$\p(a_0, \dots , a_n)$;
on a l'analogie suivante entre $U(a_0, \dots , a_n)$ et l'espace
projectif usuel:

- Le fibr\'e canonique: $K_U={\cal O}_U(-\sum a_i)$. Ceci peut se voir,
par exemple, en regardant la ramification de $\pi$: chaque coordonn\'ee
$x_i$ s'annule le long d'une composante de ramification de multiplicit\'e
$a_i-1$; on conclut en appliquant la formule de Hurwitz \`a $\pi$ au-dessus
de $U$. 

- La suite d'Euler: 

$$0\rightarrow {\cal O}_U\rightarrow \oplus {\cal O}_U(a_i)\rightarrow T_U
\rightarrow 0.$$

Elle provient du fait que, tout comme l'espace tangent \`a l'espace projectif
ordinaire est engendr\'e par les $\partial/\partial x_i$ 
(i.e. les images directs 
des $\partial/\partial x_i$ usuels sur $\C^{n+1}$!), avec la seule relation
$\sum x_i\partial/\partial x_i=0$, l'espace tangent \`a notre espace
projectif \`a poids est engendr\'e par les $\partial/\partial x_i$
avec la seule relation $\sum d_ix_i\partial/\partial x_i=0,$
o\`u $d_i=a_0\dots a_{i-1}a_{i+1}\dots a_n$.

Les vari\'et\'es de Fano de petit degr\'e (celles dont
le g\'en\'erateur $H$ n'est pas tr\`es ample) de la liste
ci-dessus peuvent \^etre d\'ecrites en termes des espaces projectifs
\`a poids. En effet:

- Un rev\^etement double de $\p^n$ ramifi\'e le long d'une 
hypersurface de degr\'e $2k$ d'\'equation $f(x_0,\dots, x_n)=0$, est donn\'e 
par l'\'equation de degr\'e $k$
$$y^2=f(x_0,\dots, x_n)$$ dans $\p(1,1,\dots ,1,k)$, avec $x_i$ les
coordonn\'ees de poids 1 et $y$ la coordonn\'ee de poids $k$; 
c'est donc une hypersurface de degr\'e $2k$ dans $\p(1,1,\dots ,1,k)$.

- Le c\^one sur la surface de Veronese est $\p(1,1,1,2)$, plong\'e dans $\p^6$ par 
$|{\cal O}_{\p(1,1,1,2)}(2)|$; une section
hypercubique de ce c\^one est une hypersurface de degr\'e 6 
dans cet espace, disons d'\'equation $g(x_0,x_1,x_2, y)=0$. 
Le rev\^etement double du c\^one ramifi\'e le long de cette
hypersurface, est donn\'e par l'\'equation de degr\'e $6$
$$z^2=g(x_0,x_1,x_2, y)$$ dans $\p(1,1,1,2,3)$ ($deg(x_i)=1$, $deg(y)=2$,
$deg(z)=3$).

- Un rev\^etement double de la quadrique dans $\p^4$ ramifi\'e le long
d'une section hyperquartique, est une intersection compl\`ete de type
(2,4) dans $\p(1,1,1,1,1,2)$ (de la quadrique dans $\p(1,1,1,1,1,2)$
passant par le point singulier et la quartique dans $\p(1,1,1,1,1,2)$
ne passant pas par le point singulier).

\

{\it 1.3. Familles de droites}

\

Soit $X$ une vari\'et\'e de Fano, $Pic(X)=\Z[H]$, $H$ ample.
On appellera une {\it droite} sur $X$ une courbe rationnelle 
$C\subset X$, telle que
$C\cdot H=1$. Si $H$ est tr\`es ample, c'est, bien s\^ur, une vraie
droite. 
On va beaucoup utiliser le fait fondamental suivant, inconnu en dimension
sup\'erieure:

\

{\bf Fait:} Si $dim(X)=3$, $X$ contient des droites.

\

Ceci pourrait se voir de la classification, mais en r\'ealit\'e la
classification elle-m\^eme depend de ce fait! Le cas r\'eellement
difficile \`a d\'emontrer est celui d'indice 1, $H^3\geq 6$, et ceci
est un th\'eor\`eme de Shokurov (\cite{S}).

Dans cette discussion de familles de droites sur $X$, on 
regardera surtout le cas o\`u $H$ est tr\`es ample (nous n'aurons
ensuite besoin que de ce cas), pour \'eviter certaines
complications. En effet, sinon, il n'est m\^eme pas clair
apriori que les ``droites'' sont lisses; bien que souvent
c'est clair de la classification. Par exemple, les ``droites'' sur $A_2$
proviennent des bitangentes \`a la quartique: l'image r\'eciproque d'une
telle bitangente sur $A_2$ est la r\'eunion de deux $\p^1$ qui s'intersectent
en deux points (en g\'en\'eral), et chaqu'une est une ``droite''. On va donc supposer
que $H$ est tr\`es ample (et consid\'erer $X$ plong\'e dans $\p^N$ par $|H|$), 
si le contraire n'est pas sp\'ecifi\'e.

\

{\it Le cas d'indice 2:} Soit $l\subset X$ une
droite. Par la formule d'adjonction $K_l=K_X|_l\otimes det(N_{l,X})$,
on a: $det(N_{l,X})=0$. Donc $N_{l,X}={\cal O}_{\p^1}(a)\oplus
{\cal O}_{\p^1}(-a)$ pour un certain $a\geq 0$. On a les \'enonc\'es suivants:

1) $a$ est 0 ou 1. En effet, on a $l\subset X\subset \p^N$ et 
$$N_{l,X}\subset N_{l,\p^N}={\cal O}_{\p^1}(1)\oplus\dots\oplus {\cal O}_{\p^1}(1).$$
Donc $a\leq 1$.

{\it Remarque:} si $H$ n'est pas tr\`es ample, il y a d'autres valeurs de $a$ possible, mais
``pas beaucoup'': on peut par exemple v\'erifier que sur $A_2$, il y a, 
en plus, un nombre fini de droites du fibr\'e normale
${\cal O}_{\p^1}(2)\oplus {\cal O}_{\p^1}(-2)$, qui proviennent
des droites contenues dans la quartique.

{\it Cons\'equence:} La th\'eorie de d\'eformations dit alors que 
le sch\'ema de Hilbert param\'etrisant les droites est une surface lisse:
$h^0(l, N_{l,X})=2$, $h^1(l, N_{l,X})=0$ .

2) Le fibr\'e normal de la droite g\'en\'erique
est ${\cal O}_{\p^1}\oplus {\cal O}_{\p^1}$: ceci est \'equivalent
au fait que les droites couvrent $X$. Notre famille de droites \'etant \`a deux 
param\`etres, il est alors clair que les droites ne peuvent
pas couvrir un diviseur seulement, car la seule surface ayant une
famille de droites \`a deux param\`etres, est un plan, et un plan ne
peut pas \^etre un diviseur sur $X$ du fait que $Pic(X)=\Z$.

3) Il existe une famille \`a un param\`etre de droites du fibr\'e
normal ${\cal O}_{\p^1}(1)\oplus {\cal O}_{\p^1}(-1)$. Plus pr\'ecisement,
soit ${\cal H}$ le sch\'ema de Hilbert des droites sur $X$, 
$q:F\rightarrow {\cal H}$ la famille universelle des droites sur $X$
et $p:F\rightarrow X$ le morphisme naturel (d'\'evaluation).
$X$ est simplement connexe et $Pic(X)=\Z$ (tandis que $F$ est fibr\'e en droites), 
on en d\'eduit que l'application $p$ a un lieu critique de
codimension un. Consid\'erons $dp$ le long d'une droite $l$; il induit un morphisme
$$dp: N_{l,F}=T_{q(l)}{\cal H}\otimes {\cal O}_l\rightarrow N_{l,X},$$ et par
la th\'eorie de d\'eformations, ceci est pr\'ecisement l'application
d'\'evaluation $$H^0(l, N_{l,X}) \otimes {\cal O}_l \rightarrow N_{l,X}.$$
On voit donc que les droites dont le fibr\'e normal est trivial, ne peuvent
pas \^etre dans le lieu critique (car $dp$ est alors un isomorphisme); ce lieu
est donc form\'e des droites du fibr\'e normal ${\cal O}_{\p^1}(-1)\oplus {\cal O}_{\p^1}(-1)$.
Sur $X$, ils couvrent un diviseur $D$ caracteris\'e par la propri\'et\'e suivante:
{\it par un point g\'en\'erique de $D$ passe moins de droites que par un point
g\'en\'erique de $X$}.

{\it Remarque:} La {\it simple-connexit\'e} d'un $X$ de Fano d'indice 2 avec $H_X$ tr\`es
ample est claire de la classification (en appliquant le th\'eor\`eme de
section hyperplane de Lefschetz); sans classification, elle r\'esulte des \'enonc\'es classiques: en effet, une section hyperplane lisse 
de $X$
est une surface de del Pezzo (c'est un nom sp\'ecial pour une vari\'et\'e
de Fano de dimension 2), dont la rationalit\'e est classiquement connue, en
particulier, une telle surface est simplement connexe. L'application
du th\'eor\`eme de section hyperplane de Lefschetz finit l'argument. 
En fait, il est bien vrai
que {\it toute} vari\'et\'e de Fano est simplement connexe; c'est un 
r\'esultat non-trivial. Nous
nous r\'ef\'erons \`a \cite{D} pour une preuve algebro-g\'eom\'etrique et des
r\'ef\'erences \`a une preuve par la g\'eom\'etrie diff\'erentielle.

\

{\it Le cas d'indice 1:} Par adjonction, $det(N_{l,X})={\cal O}_{\p^1}(-1)$.
On a bien l'analogue de la propri\'et\'e 1:

1) Il y a deux possibilit\'es: soit $N_{l,X}=  
{\cal O}_{\p^1}(-1)\oplus {\cal O}_{\p^1}$, soit $N_{l,X}=  
{\cal O}_{\p^1}(-2)\oplus {\cal O}_{\p^1}(1)$
(m\^eme argument que pour l'indice deux).

On en d\'eduit par la th\'eorie de d\'eformations, que le sch\'ema de Hilbert
des droites est une courbe (les droites ne peuvent pas couvrir $X$ au cause
de facteurs n\'egatifs dans le fibr\'e normal), que les droites
du fibr\'e normal ${\cal O}_{\p^1}(-1)\oplus {\cal O}_{\p^1}$ correspondent
a ses points lisses et les droites du fibr\'e normal 
${\cal O}_{\p^1}(-2)\oplus {\cal O}_{\p^1}(1)$, aux points singuliers.

On {\it n'a pas} d'analogue de la propri\'et\'e 2: il se peut que 
toute droite a pour fibr\'e normal ${\cal O}_{\p^1}(-2)\oplus {\cal O}_{\p^1}(1)$
(le sch\'ema de Hilbert est non-r\'eduite). Mais on peut \^etre
un peu plus pr\'ecis ici:

2) Soit $F\rightarrow {\cal H}_{red}$ (une composante irr\'eductible de) 
la famille universelle des droites sur $X$.
Si une droite g\'en\'erique a pour fibr\'e normal ${\cal O}_{\p^1}(-1)\oplus {\cal O}_{\p^1}$,
alors le morphisme $p:F\rightarrow X$ est une immersion le long de la droite  
g\'en\'erique. Si le fibr\'e normal de toute droite est 
${\cal O}_{\p^1}(-2)\oplus {\cal O}_{\p^1}(1)$, alors $p$ a un lieu critique coupant
la droite g\'en\'erique en un point, et l'image de $F$ est soit un c\^one, soit
une surface tangente \`a une courbe.

Pour la preuve on fait le m\^eme argument qu'avant: on considere $dp$ le long $l$
qui est ou bien $dp:{\cal O}_l\rightarrow {\cal O}_{\p^1}(-1)\oplus {\cal O}_{\p^1}$,
ou bien $dp:{\cal O}_l\rightarrow {\cal O}_{\p^1}(-2)\oplus {\cal O}_{\p^1}(1)$.

3) Quelques commentaires: d'abord, une $V_k$ avec $k\geq 6$ ne contient pas de c\^ones.
En effet, un tel c\^one de sommet $x$ doit \^etre contenu dans $T_xX$; mais $X$ est une
intersection de quadriques et de cubiques, donc le degr\'e du c\^one serait au plus 3,
ce que contredit \`a $Pic(X)=\Z$. Ensuite, pour une vari\'et\'e de Fano g\'en\'erique
dans sa famille, ${\cal H}$ est bien r\'eduite (mais nous n'allons pas le demontrer ici).
Enfin, le seul exemple connu d'une vari\'et\'e de Fano d'indice 1 contenant une surface
tangente \`a une courbe, est la vari\'et\'e de type $V_{22}$ construite par Mukai
et Umemura dans \cite{MU}; la surface est une section hyperplane de cette vari\'et\'e.
Iliev et Schuhmann (\cite{IS}) ont d\'emontr\'e que les autres vari\'et\'es de Fano d'indice 1
ne peuvent pas contenir une surface tangente \`a une courbe en tant qu'une
section hyperplane, c.\`a.d. le sch\'ema de Hilbert ne peut pas \^etre
``maximalement non-r\'eduite'', sauf pour la vari\'et\'e de Mukai-Umemura.

\

Dans ce qui suit, on va appeller  $(a,b)$-droites les droites dont le fibr\'e normal
est ${\cal O}_{\p^1}(a)\oplus {\cal O}_{\p^1}(b)$.

\

{\it 1.4. Exemples de calculs}

\

Donnons quelques exemples de calculs des invariants de familles de droites
sur certaines vari\'et\'es de Fano. On ne fera pas toutes les justifications
necessaires, parfois un peu longues; ce que nous int\'eresse ici, c'est de voir
que ces invariants peuvent \^etre calcul\'es sans grande difficult\'e.

{\it Exemple 1: Droites sur une cubique dans $\p^4$.} 
Identifions d'abord la surface de droites
${\cal H}\subset G(1,4)$. La famille universelle des droites dans $\p^4$ est
la projectivisation $\p(U)$ du fibr\'e universel $U$ sur $G(1,4)$. Soient
$p$ resp. $q$ les projectons de la famille universelle vers $\p^4$ resp. $G(1,4)$.
Les droites
contenues dans la cubique donn\'ee forment le lieu des z\'eros de la section 
correspondante de $q_*p^*{\cal O}_{\p^4}(3)=S^3U^*$. En tant que cycle sur
$G(1,4)$, ${\cal H}$ est donc num\'eriquement \'equivalente \`a $c_4(S^3U^*)$. 
Ce dernier se calcule par un proc\'ed\'e
standart \`a partir de classes de Chern $c_i$ de $U^*$: 
$$c_4(S^3U^*)=9c_2(2c_1^2+c_2).$$
Notons que $$c_1(U^*)=\sigma_1=\{droites\ coupant\ un\ plan\ donn\acute{e} \},$$
$$c_2=\sigma_{1,1}=\{droites\ contenues\ dans\ un\ hyperplan\ donn\acute{e} \},$$
$$c_1^2=\sigma_{1,1}+\sigma_2,\ \sigma_2=\{droites\ coupant\ une\ droite\ donn\acute{e}e \}.$$

Pour trouver le nombre de droites passant par un point g\'en\'erique
de $X$, il suffit de calculer $1/3({\cal H}\cdot\sigma_2)$. Il est clair
que $c_2\cdot c_2\cdot \sigma_2=0$ et $c_2\cdot\sigma_2\cdot\sigma_2=1$,
on a donc 6 droites par un point g\'en\'erique.

Trouvons maintenant le degr\'e de la surface $D$ couverte sur $X$ par les
(-1,1)-droites. Soit $l$ une droite g\'en\'erique et $C_l$ une courbe des
droites coupant $l$ et diff\'erents de $l$ ($C_l$ est compl\`ete d\`es que
$l$ est un (0,0)-droite). Toutes les $C_l$ sont num\'eriquement 
\'equivalentes en tant que diviseurs sur ${\cal H}$ (elles le sont m\^eme alg\'ebriquement:
on peut d\'emontrer que ${\cal H}$ est connexe, donc, pour toutes $l,l'$, 
$C_l$ se d\'eforme en $C_{l'}$). L'intersection de $l$ avec $D$ est
l'ensemble des valeurs critiques de l'application $C_l\rightarrow l$ envoyant
$l'\in C_l$ en $l'\cap l$ (``les points par lesquels il passe moins de droites
qu'il faudrait''). Cette application est de degr\'e 5 (car il y a 6 droites par un point
g\'en\'erique). 

Par adjonction, $K_{\cal H}=K_{G(1,4)}\otimes det(S^3U^*)|_{\cal H}=
{\cal O}_{\cal H}(1)$. Evidemment, $3C_l\in |{\cal O}_{\cal H}(1)|$
(il y a des sections planes de $X$ reunions de 3 droites). On calcule
donc $degK_{C_l}$ par adjonction: $degK_{C_l}=4C_l^2$; et on remarque
que $C_l^2=5$ (6 droites par un point g\'en\'erique donc 5 droites coupant
2 droites qui s'intersectent). On en deduit que la ramification du
$C_l\rightarrow l$ est de degr\'e 30; on expecte donc que 
$D\sim 30H_X$. En tout cas, m\^eme si la ramification n'est pas simple, on aura
$D\sim kH_X$ avec $k\geq 8$.

\

{\it Exemple 2: droites sur $A_5$.} Traitons d'une autre mani\`ere la vari\'et\'e 
$A_5=G(1,4)\cap H_1\cap H_2\cap H_3$.
Les droites sur $A_5$ correspondent aux cycles de Schubert 
$\sigma_p=\p^3\subset G(1,4)$,
o\`u
$$\sigma_p=\{droites\ passant\ par\ le\ point\ p\in\p^4\},$$ tels que 
$\sigma_p\cap H_1\cap H_2\cap H_3$ est de dimension 1 (au lieu de dimension 0 
g\'en\'erique). Notre surface ${\cal H}$ est donc naturellement plong\'ee dans $\p^4$.
On peut voir ce plongement autrement: le $\p^9$ dans lequel on a plong\'e $G(1,4)$ par
les coordonn\'ees de Pl\"ucker, est $\p(\Lambda^2\C^5)$, les hyperplans
sont donc des \'el\'ements de $\p(\Lambda^2\C^{5*})$, c'est-\`a-dire, les formes
antisym\'etriques sur $\C^5$. Une telle forme $A$ est g\'en\'eriquement de rang 4, et a
donc une droite vectorielle - un point de $\p^4$! - pour noyau, le ``centre'' de $A$.
Pour $p,q\in\p^4$ $(p\neq q)$, $pq\in H_A$ \'equivaut \`a $\ ^tpAq=0$. Donc, $p$ est le centre
de $A$ si et seulement si toute droite par $p$ est dans $H_A$. Ceci implique 
que notre surface ${\cal H}$ est form\'ee de centres des hyperplans de notre reseau
$\lambda H_1+\mu H_2 +\nu H_3$, donc isomorphe \`a $\p^2$. Des calculs classiques
(voir par exemple \cite{SR}) montrent que ${\cal H}$ est une projection d'une surface de Veronese
dans $\p^5$: effectivement, le noyau de $\lambda A_1+\mu A_2+\nu A_3$ s'exprime \`a l'aide des
formes quadratiques en $\lambda, \mu, \nu$.

Notre vari\'et\'e $A_5$ est naturellement (la cl\^oture dans $G(1,4)$ de) 
l'ensemble des droites 3-secantes
de ${\cal H}$: effectivement, si une droite contient trois centres, le point correspondant de $G(1,4)$ 
est dans les trois hyperplans.
En particulier, on a 3 droites par un point g\'en\'erique de $A_5$. Le lieu o\`u il y a 
``moins de droites que d'habitude'' correspond \`a l'ensemble des ``3-secantes d\'eg\'en\'er\'ees''.
Pour identifier ce lieu, diviseur $D$ (couvert sur $A_5$ par (-1,1)-droites), reinterpr\'etons
tout en termes de {\it triangles polaires}.

Si $p_1, p_2, p_3$ sont trois points align\'es de ${\cal H}\subset \p^4$, alors
les points correspondants $p_1', p_2', p_3'$ de la surface de Veronese $v_2(\p^2)$
sont dans un plan passant par le centre $P$ de la projection. Mais l'application
de Veronese $v_2$ plonge $\p^2=\p(V)$ dans $\p^5=\p(S^2V^*)$, l'espace des coniques sur $\p^2$, et
les points de $v_2(\p^2)$ correspondent aux droites doubles. Le fait
que $p_i'$ et $P$ sont coplanaires se traduit donc par
$$Q=l_1^2+l_2^2+l_3^2,$$
o\`u les $l_i$ sont des formes lin\'eaires repr\'esentant $p_i'$ et $Q$ est une
forme quadratique repr\'esentant $P$. On peut donc interpr\'eter
$A_5$ comme la {\it vari\'et\'e des triangles polaires} d'une conique plane,
la cl\^oture de l'ensemble des triangles polaires dans $Hilb^3(\p(V^*))$.
Le diviseur $D$ est forme des ``triangles polaires d\'eg\'en\'er\'es'', dont un calcul
facile montre que ce sont
des r\'eunions d'une droite tangente \`a la conique et une droite passant par
le point de tangence. Les droites sur $A_5$ sont les ensembles des triangles
polaires incluant une droite donn\'ee. Une droite g\'en\'erique contient,
bien s\^ur, deux points correspondant aux triangles ``d\'eg\'en\'er\'es'',
donc $D$ est une section hyperquadrique de $A_5$. (En fait, un calcul simple montre que
c'est la surface tangente \`a une courbe rationnelle normale de degr\'e 6). 

\

{\it 1.5. Coniques sur les vari\'et\'es de Fano d'indice 1}

\

On donne ici, sans d\'etailler les preuves, quelques
propri\'et\'es de familles de coniques sur un $X$ d'indice 1 plong\'e
anticanoniquement dans $\p^N$. 

1) L'existence: elle peut se d\'eduire de l'existence de droites.
Plus pr\'ecisement: on demontre que si $D$ est un diviseur sur $X$
couvert par des droites dont la g\'en\'erique est (0,-1), alors
ces droites se coupent (une droite g\'en\'erique intersecte
$d+1$ droites, o\`u $d$ est tel que $D\sim dH_X$). Si l'on a des (0,-1)-droites, on
a donc des coniques r\'eductibles. Sinon, on montre que sur une
(1,-2)-droite il y a une structure double localement intersection
complete de genre arithmetique 0 (``conique''). Ensuite, par la
th\'eorie de d\'eformations, on voit que la dimension de la famille
de coniques est au moins deux (en particulier il y a des coniques lisses).

2) Le fibr\'e normal: pour une conique lisse, c'est 
${\cal O}_{\p^1}(a)\oplus {\cal O}_{\p^1}(-a)$ avec $a=0,1,2$ ou $4$.
Ceci s'ensuit, comme avant, de 
$$N_{C,X}\subset N_{C, \p^N}={\cal O}_{\p^1}(4)\oplus {\cal O}_{\p^1}(2)^{\oplus N-2}.$$
En plus, on en d\'eduit que si $a=4$, $X$ est une quartique et le plan engendr\'e par $C$ est
tangent \`a $X$ le long $C$. Pour une conique r\'eductible $l_1\cup l_2$,
la r\'estriction du fibr\'e normal sur $l_i$ est ${\cal O}_{\p^1}(a)\oplus {\cal O}_{\p^1}(-a)$
avec $a=0,1, 2$. Si $a=2$ pour les deux $i$, alors $X$ est une quartique et le plan contenant
$l_1\cup l_2$ est tangent \`a $X$ le long $l_1\cup l_2$.

3) Les coniques couvrent $X$, une conique g\'en\'erique est donc (0,0).

\

{\bf Partie II: Quelques m\'ethodes pour obtenir des bornes pour le degr\'e
de morphismes} 

\

Le th\'eor\`eme suivant r\'esulte de \cite{A2}, \cite{ARV}:

\

{\bf Th\'eor\`eme:} {\it Soient $X, Y$ deux vari\'et\'es de dimension 3, $b_2(X)=b_2(Y)=1$.
Si $Y\not\cong \p^3$, alors le degr\'e d'un morphisme $f:X\rightarrow Y$ est born\'e
par des invariants discrets de $X,Y$.}

\

Ce th\'eor\`eme est le plus difficile \`a d\'emontrer dans le cas o\`u $Y$ est Fano;
la preuve depend de la classification et de la g\'eom\'etrie des droites sur $Y$. 

Dans cette partie, $Y$ sera une vari\'et\'e de Fano ($\not\cong \p^3$) de dimension 3 avec le groupe de
Picard cyclique; $X$ sera une vari\'et\'e quelconque avec le groupe de Picard cyclique;
$f:X\rightarrow Y$ sera un morphisme. On discutera des divers moyens pour montrer
que $deg(f)$ est born\'e.

\

{\it 2.1. Classes de Chern}

\
 
Supposons pour un instant que, au lieu de supposer que $-K_Y$ est ample,
on suppose que $K_Y$ est ample. Dans ce cas, le fait que $deg(f)$ est born\'e
est compl\'etement trivial en toute dimension 
et r\'esulte de ce qu'une forme diff\'erentielle,
relev\'ee de $Y$ \`a $X$, peut acquerir des z\'eros mais pas en perdre, ni acquerir des
p\^oles. 
En effet, \'ecrivons $K_X=f^*K_Y+R$; $R$ est un diviseur effectif. Si $n=dim(X)=dim(Y)$,
\'ecrivons $K_X^n=(f^*K_Y+R)^n$; comme $f^*K_Y$ est ample et $R$ est effectif, ceci
donne $K_X^n\geq (f^*K_Y)^n = deg(f)K_Y^n$; voici une borne.

Cet argument trivial ne marchera pas si $K_Y$ est nul, ou, a fortiori, n\'egatif;
mais on peut esperer d'obtenir quelque r\'esultat en considerant non pas la classe canonique
qui est le d\'eterminant du fibr\'e cotangent, mais des classes de Chern sup\'erieures
de ce m\^eme fibr\'e.

\

{\bf Lemme 2.1.1} (\cite{ARV}) {\it Soient $X$, $Y$ lisses projectives de dimension $n$, $f:X\rightarrow Y$ 
un morphisme fini. Soit $L$ un fibr\'e en droites sur $Y$ tel que $\Omega_Y(L)$ est engendr\'e
par ses sections globales. Alors, pour $s\in H^0(Y,\Omega_Y(L))$ g\'en\'erique, les z\'eros de la section
induite $\tilde{s}\in H^0(X,\Omega_X(f^*L))$ sont isol\'es.}

\

{\it D\'emonstration}: Pour simplicit\'e, on va traiter le cas $n=2$, le cas g\'en\'eral \'etant
analogue. D'abord, les z\'eros de $s$ sont bien isol\'es parce que $\Omega_Y(L)$ est engendr\'e
par ses sections. Puis, $\tilde{s}$ a des z\'eros aux points de l'image r\'eciproque des z\'eros de $s$
(z\'eros de $f^*s$), plus , peut-\^etre, sur la ramification: le passage de
$f^*s$ \`a  $\tilde{s}$ est donn\'e par la matrice jacobienne de $f$. Si les z\'eros de $\tilde{s}$ ne
sont pas isol\'es, il y a donc une composante $R^0$ de la ramification le long de laquelle une
$\tilde{s}$ g\'en\'erique (et, par cons\'equent, toute $\tilde{s}$) s'annule. 

Consid\'erons $\phi: f^*(\Omega_Y(L))\rightarrow \Omega_X(f^*L)$, l'application donn\'ee localement
par la matrice jacobienne. Au point g\'en\'erique de $R^0$, le rang de  $\phi$ est 1 (car notre
morphisme est fini). Donc, le fait que   $\tilde{s}$ s'annule le long de $R^0$ signifie que
la restriction de $f^*s$ sur $R^0$ est dans un sous-faisceau de rang un de $f^*(\Omega_Y(L))|_{R^0}$,
``le noyau de la jacobienne''. Mais ceci est impossible: les $f^*s$ engendrent $f^*(\Omega_Y(L))$,
et leurs restrictions sur $R^0$ engendrent $f^*(\Omega_Y(L))|_{R^0}$.

Le cas g\'en\'eral se fait de la m\^eme mani\`ere, en observant que $$dim\{x\in X|rg_x(f)\leq i\} \leq i.$$

\

{\bf Corollaire 2.1.2:} 
$$deg(f)c_n(\Omega_Y(L))\leq c_n(\Omega_X(f^*L)).$$

\

{\it Preuve:} Choisissons $s$ g\'en\'erique, aux z\'eros isol\'es qui sont en dehors du lieu des valeurs
critiques de $f$ (ceci est possible car $\Omega_Y(L)$ est engendr\'e). Alors $deg(f)c_n(\Omega_Y(L))$
est le nombre des z\'eros de $f^*s$ et $c_n(\Omega_X(f^*L))$ est le nombre des z\'eros de $\tilde{s}$.
Clairement, $\tilde{s}$ a un z\'ero partout o\`u $f^*s$ en a.

\

Re\'ecrivons cette in\'egalit\'e dans le cas o\`u $n=3$, $Pic(X)=Pic(Y)=\Z$.
Soit donc $l$ tel que $L\sim lH_Y$ et $m$ tel que $f^*H_Y=mH_X$; le degr\'e
de $f$ est alors $m^3H_X^3/H_Y^3$. On a:

$$deg(f)(c_3(\Omega_Y)+c_2(\Omega_Y)lH_Y+c_1(\Omega_Y)l^2H_Y^2+l^3H_Y^3)\leq$$
$$\leq c_3(\Omega_X)+c_2(\Omega_X)lmH_X+c_1(\Omega_X)l^2m^2H_X^2+l^3m^3H_X^3$$

Les derniers termes s'annulent; apr\`es cette annulation, la partie gauche
croit comme $m^3$ et la partie droite comme $m^2$. Donc, si 
$$E(Y,l)=c_3(\Omega_Y)+c_2(\Omega_Y)lH_Y+c_1(\Omega_Y)l^2H_Y^2 > 0,$$
on obtient que $m$ et donc aussi $deg(f)$ est born\'e.

Pour maximiser cette derni\`ere expression, on a bien s\^ur l'int\'er\^et de chercher
le plus petit $l$ tel que $\Omega_Y(L)$ est engendr\'e, puisque le terme en $l^2$
est n\'egatif. Si $H_Y$ est tr\`es ample, soit $Y$ plong\'e dans $\p^N$ par $H$;
alors $\Omega_Y(L)$ est quotient de $\Omega_{\p^N}(l)|_Y$. On a 
$$\Omega_{\p^N}\cong \Lambda^{N-1}T_{\p^N}\otimes K_{\p^N} = 
\Lambda^{N-1}T_{\p^N}(-N-1)$$
et $\Lambda^{N-1}T_{\p^N}(-N+1)=\Lambda^{N-1}(T_{\p^N}(-1))$ est engendr\'e;
donc on peut prendre $l=2$.

Si $H_Y$ n'est pas tr\`es ample, on peut, bien s\^ur, prendre $l=2k$ o\`u $k$ est tel
que $kH_Y$ est tr\`es ample, mais ceci ne suffit pas dans certains cas int\'eressants.
C'est ici qu'on aura recours aux espaces projectifs \`a poids.

\

{\it Exemples:}

1) Quartique: pour $Y$ une quartique, $b_3(Y)=60$ donc $c_3(\Omega_Y)=\chi(Y)=56$;
$c_2(\Omega_Y)\dot H= c_2c_1 =24$. Notre expression $E(Y,2)$ est donc
$56+48-16>0$.

2) $Y=A_4$: $b_3=4$, $c_3(\Omega_Y)=0$, $c_2(\Omega_Y)\dot 2H= c_2c_1 =24$;
$E(Y,2)=24-32<0$, ceci ne marche donc pas!

3)Regardons un exemple o\`u $H$ n'est pas tr\`es ample, disons $Y=A_2$.
$\Omega_Y(L)$ est un quotient de $\Omega_U(l)|_Y$, o\`u $U$ est la partie lisse
de $\p(1,1,1,1,2)$. On a:
$$\Omega_U\cong \Lambda^3T_U\otimes det(\Omega_U)= \Lambda^3T_U(-6)$$
(la classe anticanonique est la somme des poids). Par la suite d'Euler,
$T_U$ est un quotient de ${\cal O}(1)^{\oplus 3}\oplus {\cal O}(2)|_U$.
En plus ${\cal O}(1)|_U$ est engendr\'e par ses sections globales: le seul
point de base de l'espace de ses sections sur $\p(1,1,1,1,2)$ est le point singulier $(0:0:0:0:1)$.
Donc $\Lambda^3T_U(-3)$ est engendr\'e, et alors on peut prendre $l=3$.

Calculons: $b_3=20$, $c_3(\Omega_Y)=16$; $E(Y,3)=16+24-36>0$.

Notons que prendre $l=4$ ne suffit pas ($E(Y,4)$ sera nulle!), donc voici l'utilit\'e
des espaces projectifs \`a poids.

\

En fait les calculs montrent que pour toutes nos vari\'et\'es dont $H$ n'est pas tr\`es
ample, la formule marchera. En plus, elle marche pour la cubique et pour $V_k, k\leq 12$,
mais ne marche pas pour les vari\'et\'es qui restent.

\

{\it Remarque:} Cette formule est efficace dans d'autres situations aussi.
Par exemple, on peut montrer que si le degr\'e d'un morphisme $f:X\rightarrow Y$,
o\`u $X, Y$ sont de dimension quelconque avec Picard $\Z$ et $K_Y=0$, n'est pas born\'e,
alors $Y$ est un quotient d'un tore par un groupe fini. On peut aussi montrer que
les hypersurfaces lisses de degr\'e $>2$ dans $\p^{>2}$ n'ont pas d'endomorphismes de degr\'e $>1$ 
(Beauville), ainsi que trouver des contraintes pour le degr\'e d'un morphisme entre
hypersurfaces.

\

{\it 2.2. Droites aux facteurs n\'egatifs dans le fibr\'e normal}

\

Cette m\'ethode marche s'il y a ``beaucoup'' de courbes dont le fibr\'e
normal a une composante (facteur direct) n\'egative, par exemple, de droites sur $Y$ d'indice 1
ou de (-1,1)-droites sur $Y$ d'indice 2.

Pr\'ecisons que veut dire ``beaucoup''. Soit $f:X\rightarrow Y$ un morphisme fini,
$Pic(X)=Pic(Y)=\Z$,
$S\subset Y$ la surface couverte par nos courbes. ``Beaucoup'' signifie que,
d\`es que $deg(f)>> 0$, on peut garantir que $f^{-1}(S)$ n'est pas contenue
dans la ramification, i.e. que l'image r\'eciproque (sch\'ematique) d'une courbe 
g\'en\'erique
de notre famille a une composante r\'eduite.

\

{\bf Lemme 2.2.1:} {\it Cette condition est bien remplie:
a) si $Y$ est d'indice 1, $S$ est la surface couverte par les droites,
et $S\sim kH_Y$ avec $k\geq 3$;
b)  si $Y$ est d'indice 2, $S$ est la surface couverte par les (-1,1)-droites,
et $S\sim kH_Y$ avec $k\geq 5$.

Si l'on suppose que $X$ lui aussi est de Fano, on peut remplacer $k\geq 3$
par $k\geq 2$ et $k\geq 5$ par $k\geq 4$, et sous ces conditions $f^{-1}(S)$ n'est jamais
contenue dans le lieu critique (m\^eme sans restriction $deg(f)>> 0$).}

\

{\it Preuve:} Consid\'erons par exemple le cas d'indice 1. Si $S\sim kH_Y$ avec $k\geq3$,
et $m$ est tel que $f^*H_Y=mH_X$,
alors $f^*S\sim mkH_X$. Si $f^{-1}S=Supp(f^*S)$ est contenu dans le lieu de ramification,
alors le diviseur de ramification $R\sim lmH_X$ avec $l\geq 3/2$: en effet, une composante
entrant en $f^*S$ avec la multiplicit\'e $k$, entre en $R$ avec la multiplicit\'e $k-1$,
et $k$ est toujours $\geq 2$. On a $K_X=-mH_X+R$, d'o\`u: soit $m$ est petit, soit 
$f^{-1}S$ n'est pas contenu dans le lieu critique. 

Si $X$ est de Fano, donc $K_X$ est n\'egatif,
alors $R\sim lmH_X$ avec $l<1$; donc d\'ej\`a l'image reciproque de $S\sim 2H_Y$ ne peut
pas \^etre contenue dans le lieu critique: on aurait alors, par le m\^eme raisonnement 
qu'avant, $l\geq 1$. 

Le cas d'indice deux est analogue. 

\

{\bf Proposition 2.2.2} {\it Le fait qu'il y a ``beaucoup'' de courbes comme ci-dessus
implique que $deg(f)$ est born\'e.}

\ 

{\it Preuve:} Soit, par exemple, $X$ d'indice un avec ``beaucoup''
de droites. Il existe alors une droite $l\subset X$, telle que $C=f^*l$ (l'image r\'eciproque
sch\'ematique de $l$) a une composante irr\'eductible r\'eduite $D$.
L'inclusion de faisceaux d'id\'eaux $I_C\rightarrow I_D$ induit le morphisme
de faisceaux conormaux $(I_C/{I_C}^2)\otimes {\cal O}_D\rightarrow I_D/{I_D}^2$, qui
est un isomorphisme partout o\`u $D$ n'intersecte pas d'autres composantes de $C$.
On a aussi un morphisme $T_X|_D\rightarrow (I_D/{I_D}^2)^*$, surjectif au point lisse de $D$.
Ensemble, ceci donne une surjection g\'en\'erique
$$T_X|_D\rightarrow ((I_C/{I_C}^2)|_D)^* = {\cal O}_D(-m)\oplus {\cal O}_D\ ou\ 
{\cal O}_D(-2m)\oplus {\cal O}_D(m),$$
selon ce que $l$ est (0,-1)- ou (1,-2)-droite.
Il est donc clair que $m$ est born\'e par tout $j$ tel que $T_X(j)$ est engendr\'e
par ses sections globales; et il existe des moyens d'exprimer un tel $j$ en termes
des invariants discr\`ets de $X$ (par exemple, si $X$ est une hypersurface de degr\'e
$d$, on peut prendre $j=d-2$. En g\'en\'eral, le probl\`eme le plus difficile est de savoir
quel multiple de $H_X$ est tr\`es ample; il existe de nombreuses r\'ef\'erences sur le
sujet, par exemple, \cite{Dem}).

Le cas d'indice deux est encore compl\'etement analogue.

\

{\it Commentaire:} On veut, bien s\^ur, savoir sur quelles vari\'et\'es de Fano
il y a ``beaucoup'' de courbes sp\'eciales. Il apparait que pour l'indice deux,
seulement $A_5$ ne satisfait pas \`a la condition du lemme: $S\sim 2H$ (comme
on l'a calcul\'e dans la premi\`ere partie). Pour l'indice 1, c'est plus 
d\'elicat: les calculs comme ceux qu'on a fait dans la premi\`ere partie,
donneraient $deg(S)$ ``en comptant avec les multiplicit\'es'', et ce sera
le vrai $deg(S)$ seulement dans le cas o\`u le schema de Hilbert des droites
est r\'eduite, i.e. la droite
g\'en\'erique est (0,-1). On peut calculer que ce $deg(S)$ ``virtuel'' est
assez grand pour toute famille de vari\'et\'es sauf $V_{22}$ o\`u
$S\sim 2H$ ``virtuellement'' (et aussi r\'eellement pour toute $V_{22}$ sauf
celle de Mukai-Umemura o\`u $S\sim H$): voir \cite{I}, \cite{Mar}. Donc, la borne g\'en\'erale marche
pour un $V_k$ g\'en\'erique dans sa famille, avec $k\neq 22$.

Quant au cas o\`u $X$ lui aussi est Fano, les conditions du lemme sont 
remplies pour tout $V_k$ sauf celui de Mukai-Umemura (par un r\'esultat
d'Iliev et Schuhmann d\'ej\`a cit\'e).

\

Regardons plus en d\'etail le cas o\`u $X$ est de Fano. Dans ce cas, non
seulement les conditions du lemme peuvent \^etre affaiblies, mais on peut aussi
identifier pr\'ecisement les composantes de $f^{-1}l$.

\

{\bf Lemme 2.2.3} {\it Soit $Y$ une vari\'et\'e lisse admettant une 
courbe rationnelle $C_0$
lisse au fibr\'e normal trivial, et $f:X\rightarrow Y$ un morphisme fini
(avec $X$ lisse). Soit $C$ une d\'eformation g\'en\'erique de $C_0$;
alors $f^{-1}C$ est lisse (pas forc\'ement connexe), et son fibr\'e normal
est trivial.}

\

{\it Preuve:} Les d\'eformations de $C_0$ sont parametr\'ees par ${\cal H}$,  
$dim{\cal H}=dimX-1$; elles couvrent $Y$,
et, en considerant la diff\'erentielle $dp$ du morphisme g\'en\'eriquement
fini $p:F\rightarrow Y$ (o\`u $F$ est la famille universelle, normalis\'ee
si necessaire), on voit que la d\'eformation g\'en\'erique $C$(dont le fibr\'e
normal est toujours trivial) est transverse au lieu de valeurs critiques $B$ de $f$; 
en plus, $C$ intersecte
$B$ au points ``g\'en\'eriques'', c.\`a.d. pour tout $y\in C\cap B$ et tout $x\in f^{-1}(y)$,
l'application $f$ au voisinage de $x$ peut \^etre \'ecrit comme
$$w_1=z_1, w_2=z_2,\dots, w_{n-1}=z_{n-1}, w_n=z_n^m$$ 
(avec $m$ dependant de $x$).   
On en d\'eduit par un calcul local que
$f^{-1}C$ est lisse. Le fibr\'e normal est l'image r\'eciproque de celui de $C$. 

\

Si $X$ est de Fano, toute courbe $C$ du fibr\'e normal trivial sur $X$ est automatiquement
rationnelle, et $K_X\cdot C=-2$ (formule d'adjonction).
On d\'eduit donc du lemme 2.2.3, pour un morphisme $f:X\rightarrow Y$ de vari\'et\'es
de Fano de dimension 3 et groupe de Picard cyclique:

- Si $X,Y$ sont d'indice 1, l'image r\'eciproque d'une conique g\'en\'erique est
une r\'eunion disjointe de coniques, donc l'image r\'eciproque d'une conique ou une droite
quelconque est une r\'eunion (pas disjointe en g\'en\'eral) des coniques et des droites;

- Si $X$ est d'indice 1, $Y$ d'indice 2, l'image r\'eciproque d'une droite g\'en\'erique est
une r\'eunion disjointe de coniques, donc l'image r\'eciproque d'une droite
quelconque est une r\'eunion des coniques et des droites;

etcetera.

Il s'ensuit, par exemple, la proposition suivante:

\

{\bf Proposition 2.2.4} {\it Soient $X,Y$ vari\'et\'es de Fano de dimension 3,
Picard $\Z$, indice 1, telles que $H_X$ et $H_Y$ sont tr\`es amples, et $Y\neq V_{22}^s$
(la vari\'et\'e de Mukai-Umemura).
Alors tout $f:X\rightarrow Y$ est un isomorphisme.}

\

{\it Preuve:} Soit $m$ tel que $f^*H_Y=mH_X$. Il suffit de montrer que $m=1$:
en effet, il apparait de la classification que si $H_X^3>H_Y^3$ 
(ce qui rend possible $deg(f)>1$
m\^eme si $m=1$), alors $b_3(X)<b_3(Y)$ (ce qui rend impossible tout morphisme
$f:X\rightarrow Y$), ``le $b_3$ de $V_k$ d\'ecroit avec $k$''. Reprenons
le proc\'ed\'e de la proposition pr\'ec\'edente: soit $l$ une droite g\'en\'erique,
$C$ son image r\'eciproque, $D$ une composante irr\'eductible r\'eduite.
C'est une droite ou une conique. On a un morphisme qui est un isomorphisme
au point g\'en\'erique:
$$(I_D/I_D^2)^*\rightarrow (I_C/I_C^2)^*|_D={\cal O}_D\oplus {\cal O}_D(-m)\ ou\ 
{\cal O}_D(m)\oplus {\cal O}_D(-2m).$$ L'analyse de toutes les possibilit\'es
pour le fibr\'e normal de $C$ (il n'y en a que six) permet de conclure que
$m$ ne peut pas \^etre diff\'erent de 1.

\

On laisse au lecteur de formuler et de d\'emontrer une proposition analogue pour $ind(X)=2$, $ind(Y)=1$,
et les deux autres cas. Malhereusement, on doit toujours exclure les cas
$Y=A_5$ et $Y=V_{22}^s$, et, en plus, on doit supposer que $H_X$ et $H_Y$ sont tr\`es amples
(sinon il y a d'autres possibilit\'es pour les fibr\'es normaux, ce qui complique
l'histoire). Les cas o\`u $H_X$ ou $H_Y$ ne sont pas tr\`es amples, sont, par ailleurs,
souvent traitables par la m\'ethode des classes de Chern. J'ai elabor\'e le cas $Y=V_5$
dans un article recent; je suis convaincue que tout 
morphisme non-trivial
entre les vari\'et\'es de Fano de dimension 3 , Picard $\Z$, indice 1 ou 2, est
un rev\^etement double ``indice 1$\rightarrow$ indice 2'', ramifi\'e le long
d'un diviseur canonique, et il ne reste plus que quelques petits d\'etails (un peu ennuyeux...) 
\`a v\'erifier pour le d\'emontrer tout \`a fait. 

\

{\it 2.3. Scindage d'une suite des fibr\'es normaux}

\

On consid\`ere ici $X$ avec $H_X$ tr\`es ample, plong\'ee dans $\p^N$ par $H_X$.

Soit $C$ une courbe et $S$ une surface lisse, telles que $C\subset S\subset X$.
Considerons une suite exacte des fibr\'es normaux:
$$0\rightarrow N_{C,S}\rightarrow N_{C,X}\rightarrow N_{S,X}|_C\rightarrow 0.$$
Si $C$ est une intersection compl\`ete $S\cap S'$, $S'\subset X$, alors
$N_{C,X}\cong N_{C,S}\oplus N_{C,S'}$, $N_{C,S'}=N_{S,X}|_C$ et notre suite est
scind\'ee. Sinon, le scindage de la suite est quelque chose de tr\`es particulier.

En fait, m\^eme la surjectivit\'e de l'application sur les sections
globales
$$\alpha: H^0(C,N_{C,X}) \rightarrow H^0 (C, N_{S,X}|_C)$$
est tr\`es particuli\`ere. Car alors l'image de $\alpha$ contient l'image
de $\beta: H^0 (S, N_{S,X})\rightarrow H^0 (C, N_{S,X}|_C),$ ce qui se traduit
en termes suivants:

{\it Toute d\'eformation infinitesimale de $S$ dans $X$ contient une d\'eformation 
infinitesimale de $C$}.

Mais la m\'etaphysique g\'en\'erale dit que, dans la plupart des cas, une $S$
g\'en\'erique ``ne contient que du strict necessaire'', c'est-\`a-dire, les courbes
sur $S$ proviennent des diviseurs sur $X$.

Une illustration classique de ce principe est le

\

{\bf Th\'eor\`eme de Noether-Lefschetz:} {\it Sur une surface g\'en\'erique 
$S\subset \p^3$
de degr\'e $\geq 4$, toute courbe est une intersection compl\`ete $S\cap S'$}.

\

La diff\'erence entre cubiques et quartiques se voit d\'ej\`a au niveau de droites:
en effet, une cubique lisse contient 27 droites, tandis qu'un calcul de dimensions facile
montre qu'une quartique g\'en\'erique n'en contient pas (la codimension, dans l'espace
des quartiques, de l'ensemble des quartiques contenant une droite donn\'ee, est
$h^0(\p^1,{\cal O}_{\p^1}(4))=5$, tandis que la dimension de la famille des droites
$G(1,3)$ est 4). Le th\'eor\`eme peut se d\'emontrer, entre autres, par {\it la m\'ethode
infinitesimale:} on consid\'ere l'ouvert $U\subset\p(H^0(\p^3,{\cal O}_{\p^3}(d))$
parametrisant les surfaces lisses de degr\'e $d$, et le syst\'eme locale de la
{\it cohomologie primitive} (la cohomologie orthogonale \`a la section
hyperplane par rapport \`a la forme d'intersection) $H^2_{prim}$ la-dessus. Soit 
$S_0\subset U$ une surface lisse de degr\'e $d$, $\lambda$ une classe enti\'ere 
de type (1,1)
dans $H^2_{prim}(S_0, \Z)$, et $\Sigma_{\lambda}\in U$ le lieu o\`u 
$\lambda$ reste de type (1,1). On d\'emontre alors que pour $d\geq 4$,
$T_{S_0}\Sigma_{\lambda}$ est un sous-espace propre de $T_{S_0}U$; 
pour une courbe dont la classe dans la
cohomologie est $\lambda$, ceci signifie qu'elle ne se deforme pas avec $S_0$ en toute
direction,
m\^eme infinitesimalement. Autrement dit, on d\'eduit le th\'eor\`eme de Noether-
Lefschetz du 

\

{\bf Th\'eor\`eme de Noether-Lefschetz infinitesimal:}
{\it Soit $S\subset\p^3$ une surface lisse de degr\'e $d\geq 4$; alors les classes enti\'eres
de type (1,1) sur $S$ qui restent de type (1,1) par toute d\'eformation infinitesimale,
sont les restrictions des classes de type (1,1) sur $\p^3$.}

\

(voir par exemple \cite{CGGH} pour une d\'emonstration).

Il apparait que l'on peut remplacer $\p^3$ par une vari\'et\'e lisse projective $X$
quelconque de dimension 3. Il faut alors, bien sur, une autre hypoth\'ese
de genre ``$S$ est un diviseur assez ample sur $X$'' (au lieu de $d\geq 4$).
Ein et Lazarsfeld d\'emontrent dans \cite{EL} que la condition suivante
est suffisante:

$$S\sim 3K_X+16A,\ A\ tr\grave{e}s\ ample.$$

Les vari\'et\'es de Fano possedent
souvent des courbes dont la suite des fibr\'es normaux est scind\'ee,
on peut donc s'attendre \`a ce que le th\'eor\`eme de Noether-Lefschetz infinitesimal
ait des cons\'equences pour les morphismes vers les vari\'et\'es de Fano.
En effet, on d\'emontre, par exemple,

\

{\bf Proposition 2.3.1}{\it Soit $X$ une vari\'et\'e de dimension 3, $Pic(X)=\Z$,
et $f:X\rightarrow Q$ un morphisme sur la quadrique de dimension trois.
Alors $deg(f)$ est born\'e par des invariants discrets de $f$.}

\

{\it Preuve:} Soit $l\subset Q$ une droite, et $H\subset Q$ une section hyperplane
lisse
contenant $l$. La suite des fibr\'es normaux est scind\'ee:
$N_{l,H}={\cal O}_l$, $N_{l,X}={\cal O}_l\oplus {\cal O}_l(1)$ et 
$N_{H,X}|_l={\cal O}_l(1)$. Parce que $Q$ est homog\`ene, pour une choix
g\'en\'erique de $l$ et de $H$, les images r\'eciproques $C=f^{-1}l$ et $S=f^{-1}H$
sont lisses. La suite des fibr\'es normaux est simplement ``l'image
r\'eciproque'' de la suite des fibr\'es normaux pour $l,H,Q$; elle est donc scind\'ee.
Evidemment la classe de $C$ dans la cohomologie ne provient pas d'une classe
sur $X$: en effet, $Pic(X)=\Z$ et $C^2=0$ sur $S$. 
Ce qui implique que le th\'eor\`eme de Noether-Lefschetz infinitesimal
n'est pas vrai pour $S\subset X$. Donc $S$ ``n'est pas assez ample'', ce qui
revient \`a dire que $deg(f)$ est born\'e. La condition explicite de \cite{EL}
permet de traduire ce borne en termes des invariants de $X$.

\

Que se passe-t-il pour nos autres vari\'et\'es de Fano $Y$, qui, elles aussi,
ont souvent des droites $l\subset H\subset Y$ telles que la suite
correspondante des fibr\'es normaux est scind\'ee? 
Pour les vari\'et\'es de Fano d'indice deux avec (-1,1)-droites, et
les vari\'et\'es de Fano d'indice un avec le schema de Hilbert des droites non-reduite,
donc avec (-2,1)-droites? Malheureusement on ne peut
pas appliquer le th\'eor\`eme de Noether-Lefschetz infinitesimal tel quel,
car il y a un

{\it Probl\'em\`e: En g\'en\'eral, on ne peut pas garantir que $S=f^{-1}(H)$
est lisse.}

Effectivement, dans le cas g\'en\'eral, les $l$ sont toutes tangentes \`a une courbe $A\subset Y$;
si le rang de $f$ est 1 le long d'une composante de l'image inverse de $A$, $S$ aura une
singularit\'e \`a l'image inverse du point de tangence. N\'eanmoins, on peut toujours
demontrer un enonc\'e qui resoud notre probl\`eme, sous l'hypoth\`ese plus forte
$b_2(X)=1$: si $S$ est une surface diviseur de Cartier assez ample sur $X$, et $C\subset S\subset X$
est une courbe telle que la suite exacte des fibr\'es normaux existe et est scind\'ee,
alors $C$ est num\'eriquement \'equivalente \`a un multiple rationnel d'une classe provenant de $X$.
La condition pr\'ecise d'amplitude se traduit en termes cohomologiques: il suffit
d'avoir $h^1(S,\Omega^1_X|_S)=1$, par exemple, $h^i(X,\Omega^1_X(-S))=0$ pour $i=1,2$.
Ceci se traduit en termes d'invariants discrets de $X$, par exemple, en utilisant
l'annulation de Griffiths. 

La preuve n'a pas beaucoup en commun avec le th\'eor\`eme de Noether-Lefschetz 
infinitesimal
(voir \cite{A2} pour d\'etails).

{}

\end{document}